\documentclass[a4paper, 12pt, reqno]{amsart}
\usepackage[all]{xy}
\usepackage{amssymb,latexsym,amsmath,amsfonts}

\numberwithin{equation}{section}

\newtheorem{theorem}{Theorem}

\newtheorem{result}[theorem]{Result}
\newtheorem{proposition}[theorem]{Proposition}

\theoremstyle{remark}

\newcommand{\cO}{{\mathcal O} }
\newcommand{\wt}{\widetilde}

\begin{document}

\title[A criterion for biholomorphism]{A criterion for a degree-one
holomorphic map to be a biholomorphism}

\author[G.\ Bharali]{Gautam Bharali}

\address{Department of Mathematics, Indian Institute of Science, Bangalore 560012, India}

\email{bharali@math.iisc.ernet.in}

\author[I.\ Biswas]{Indranil Biswas}

\address{School of Mathematics, Tata Institute of Fundamental Research, 1 Homi Bhabha
Road, Mumbai 400005, India}

\email{indranil@math.tifr.res.in}

\author[G.\ Schumacher]{Georg Schumacher}

\address{Fachbereich Mathematik und Informatik,
Philipps-Universit\"at Marburg, Lahnberge, Hans-Meerwein-Stra{\ss}e, D-35032
Marburg, Germany}

\email{schumac@mathematik.uni-marburg.de}


\begin{abstract}
Let $X$ and $Y$ be compact connected complex manifolds of the same dimension with 
$b_2(X)\,=\,b_2(Y)$. We prove that any surjective holomorphic map of degree one from 
$X$ to $Y$ is a biholomorphism. A version of this was established by the first two 
authors, but under an extra assumption that $\dim H^1(X,\, {\mathcal O}_X)\,=\,\dim 
H^1(Y,\, {\mathcal O}_Y)$. We show that this condition is actually automatically 
satisfied.
\end{abstract}

\maketitle

\section{Introduction}

Let $X$ and $Y$ be compact connected complex manifolds of dimension $n$. Let
$$
f \,:\, X\,\longrightarrow\, Y
$$
be a surjective holomorphic map such that the degree of $f$ is one, meaning that
the pullback homomorphism
$$
{\mathbb Z}\,\simeq\, H^{2n}(Y,\, {\mathbb Z})\, \stackrel{f^*}{\longrightarrow}\,
H^{2n}(X,\, {\mathbb Z})\,\simeq\, {\mathbb Z}
$$
is the identity map of $\mathbb Z$. It is very natural to ask, ``Under what conditions would $f$ be a
biholomorphism?'' An answer
to this was given by \cite[Theorem~1.1]{b-b}, namely:

\begin{result}[{\cite[Theorem~1.1]{b-b}}]\label{res:b-b}
Let $X$ and $Y$ be compact connected complex manifolds of dimension $n$, and let
$f \,:\, X\,\longrightarrow\, Y$
be a surjective holomorphic map such that the degree of $f$ is one.
Assume that
\vspace{-1mm}
\begin{enumerate}
\item[$(i)$] the $\mathcal{C}^\infty$ manifolds underlying $X$ and $Y$ are diffeomorphic, and
\item[$(ii)$] $\dim H^1(X,\, {\mathcal O}_X) \,=\, \dim H^1(Y,\, {\mathcal O}_Y)$.
\end{enumerate}
\vspace{-1mm}
Then, the map $f$ is a biholomorphism.
\end{result}

In the proof of Result~\ref{res:b-b}, the condition~$(i)$ is used {\em only} in concluding
that $\dim H^2(X,\, {\mathbb R}) = \dim H^2(Y,\, {\mathbb R})$. In other words, the
proof of \cite[Theorem~1.1]{b-b} establishes that if
$$
\dim H^2(X,\, {\mathbb R})\,=\,\dim H^2(Y,\, {\mathbb R})\ \ \text{ and }\ \
\dim H^1(X,\, {\mathcal O}_X)\,=\,\dim H^1(Y,\, {\mathcal O}_Y)\, ,
$$
then\,---\,with $X$, $Y$, and $f$ as above\,---\, $f$ is a biholomorphism.

There is some cause to believe that the condition $(ii)$ in Result~\ref{res:b-b} might be
superfluous (which we shall discuss presently). It is the basis for our main theorem, which
gives a simple, purely topological, criterion for a degree-one map to be a biholomorphism:

\begin{theorem}\label{th:deg1}
Let $X$ and $Y$ be compact connected complex manifolds of dimension $n$, and let
$f \,:\, X\,\longrightarrow\, Y$
be a surjective holomorphic map of degree one. Then, $f$ is a biholomorphism if and only if
the second Betti numbers of $X$ and $Y$ coincide.
\end{theorem}

If $X$ and $Y$ were assumed to be K{\"a}hler, then Theorem~\ref{th:deg1} would follow
from Result~\ref{res:b-b}. This is because, by the Hodge decomposition, 
$\dim H^1(M,\, {\mathcal O}_M) = \frac{1}{2}\dim H^2(M,\, {\mathbb C})$ for any
compact K{\"a}hler manifold $M$. We shall show that this observation\,---\,i.e., that
condition $(ii)$ in Result~\ref{res:b-b} is automatically satisfied under the hypotheses
therein\,---\,holds
true in the general, {\em analytic} setting. In more precise terms, we have:   

\begin{proposition}\label{p:main}
Let the manifolds $X$ and $Y$ and $f \,:\, X\,\longrightarrow\, Y$ be as in
Result~\ref{res:b-b}. Then, $f$ induces an isomorphism between $H^1(X,\, {\mathcal O}_X)$ and
$H^1(Y,\, {\mathcal O}_Y)$. In particular, $\dim H^1(X,\, {\mathcal O}_X) = \dim H^1(Y,\, {\mathcal O}_Y)$.
\end{proposition}

The above proposition might be unsurprising to many. It is well known when $X$ and $Y$ are projective.
Since we could not find an explicit statement of Proposition~\ref{p:main}\,---\,and since certain supplementary
details are required in the analytic case\,---\,we provide a proof of it in Section~\ref{sec:proof}.
The {\em non-trivial} step in proving Theorem~\ref{th:deg1} uses Result~\ref{res:b-b}: given
Proposition~\ref{p:main}, our theorem follows from Result~\ref{res:b-b} and the comment above upon its proof.

\section{Proof of Proposition~\ref{p:main}}\label{sec:proof}

We begin with a general fact that we shall use several times below. For any proper holomorphic map
$F\,:\, V \,\longrightarrow\, W$ between complex manifolds, the
Leray spectral sequence gives the following exact sequence:
\begin{equation}\label{eq:leray}
0 \,\longrightarrow\, H^1(W,\, F_*\cO_V) \,\stackrel{\theta_F}\longrightarrow\, H^1(V,\, \cO_V) \,
\longrightarrow\, H^0(W,\, R^1F_*\cO_V)
\, \longrightarrow\, \cdots \, .
\end{equation}

With our assumptions on $X$, $Y$ and $f$, the map $f^{-1}$ (which is defined
outside the image in $Y$ of the set of points at which $f$ fails to be a local
biholomorphism) is holomorphic on its domain. Thus $f$ is bimeromorphic.

We note that any bimeromorphic holomorphic map of connected complex manifolds has
connected fibers, because it is biholomorphic on the complement of a thin analytic subset. In
particular, the fibers of $f$ are connected.

\noindent {\bf Claim~1.}\, {\em Let $F\,:\,V \,\longrightarrow\, W$ be a bimeromorphic
holomorphic map between compact, connected complex manifolds. The natural homomorphism
\begin{equation}\label{eq:nat}
\cO_W \,\longrightarrow\, F_*\cO_V
\end{equation}
is an isomorphism.}
\vspace{-0.5mm}

\noindent By definition, \eqref{eq:nat} is injective. In our case, it is an
isomorphism outside a closed complex analytic subset of $W$, say $\mathcal{S}$, of codimension
at least 2. So, to show that \eqref{eq:nat} is surjective, it suffices to show that given any
$w\in \mathcal{S}$, for each open connected set $U\ni w$ and each holomorphic function
$\psi$ on $F^{-1}(U)$ there is a function $H_\psi$ holomorphic on $U$ such that
$$
\psi \,=\, H_\psi\circ F \; \; \; \text{on $F^{-1}(U)$}.
$$
Since $F^{-1}$ is holomorphic on $W\!\setminus\!\mathcal{S}$, we set
$$\left.H_{\psi}\right|_{U\setminus\mathcal{S}} \,:= \,
\psi\circ (F^{-1}|_{U\setminus\mathcal{S}})\, .$$
This has a unique holomorphic extension
to $U$ by Hartogs' theorem (or more acurately: Riemann's second extension theorem),
since $\mathcal{S}$ is of codimension at least 2. As $F$ has compact, connected
fibers, this extension has the desired properties. This shows that the homomorphism in \eqref{eq:nat} is surjective.
Hence the claim.

By Claim~1, \eqref{eq:leray} yields an injective homomorphism
\begin{equation}\label{eq:le1}
\Theta_f\,:\,H^1(Y,\,\cO_Y) \,\longrightarrow\, H^1(X,\, \cO_X )\, ,
\end{equation}
which is the composition of the homomorphism $\theta_f$, as given by \eqref{eq:leray}, and the isomorphism
induced by \eqref{eq:nat}.

There is a commutative diagram of holomorphic maps
\begin{equation}\label{eq:comm}
\xymatrix{& Z \ar[d]^h \ar[dl]_g \\ X \ar[r]^f & ~Y\, , }
\end{equation}
where $h$ is a composition of successive blow-ups with smooth centers, such that the subset of $Y$
over which $h$ fails to be a local biholomorphism (i.e., the image in $Y$ of the the exceptional locus
in $Z$)
coincides with the subset of $Y$ over which $f$ fails to be a local biholomorphism. This fact
(also called ``Hironaka's Chow Lemma'')
can be deduced from Hironaka's Flattening Theorem
\cite[p.~503]{hiro}, \cite[p.~504, Corollary~1]{hiro}. We recollect briefly the argument for this.
The set $\mathcal{A}$ of
values of $f$ in $Y$ at which $f$ is not flat coincides with the set
of points over which $f$ is not locally biholomorphic. Hironaka's Flattening Theorem states that there exists
a sequence of blow-ups of $Y$ with smooth centers over $\mathcal{A}$ amounting to a
map $$h\,:\,Z \,\longrightarrow\, Y$$ such
that\,---\,with $\wt Z$ denoting the proper transform of $Y$ in $X\times_Y Z$
and ${\sf pr}_Z$ denoting the projection $X\times_Y Z \,\longrightarrow\, Z$\,---\,the map
$\wt{f} := \left.{\sf pr}\right|_{\wt Z}$
is flat. In our case this implies that $\wt{f}\,:\,\wt Z \,\longrightarrow\, Z$ is a biholomorphism.
The map $g = {\sf pr}_X\!\circ(\wt{f}\,)^{-1}$ and has the properties stated above.

The maps $h$ and $g$ above are proper modifications. Thus, all the assumptions in Claim~1 hold
true for $g\,:\,Z \,\longrightarrow\, X$. Hence, we conclude that the homomorphism
$\cO_X \,\longrightarrow\, g_*\cO_Z$ is an isomorphism. By \eqref{eq:leray} applied now to
$(V, W, F) = (Z, X, g)$, the homomorphism
\begin{equation}\label{eq:le2}
\Theta_g\,:\,H^1(X, \,\cO_X) \,\longrightarrow\, H^1(Z, \, \cO_Z),
\end{equation}
which is analogous to $\Theta_f$ above, is injective.

Similarly, the homomorphism $\cO_Y \,\longrightarrow\, h_*\cO_Z$ is an isomorphism.
Since \eqref{eq:leray}, an exact sequence, is natural, we
would be done\,---\,in view of \eqref{eq:le1}, \eqref{eq:le2} and the diagram \eqref{eq:comm}\,---\,if we show
that the homomorphism $\Theta_h\,:\,H^1(Y,\,\cO_Y) \,\longrightarrow\, H^1(Z,\,\cO_Z)$,
(given by applying \eqref{eq:leray} to $(V, W, F) = (Z, Y, h)$) is an isomorphism.

To this end, we will use the following:

\noindent {\bf Claim 2.}\, {\em For a complex manifold $W$ of dimension $n$, if
$$
\sigma\,:\,S \,\longrightarrow\, W
$$
is a blow-up with smooth center, then the direct image $R^1\sigma_* \cO_S$ vanishes.}

\noindent This claim is familiar to many.  However, since it is not so easy to point to one {\em specific}
work for a proof in the {\em analytic} case,  we indicate an argument. We first study the blow-up
$\wt \sigma\,:\, \wt S \,\longrightarrow\, \wt W$ of a point $0\in \wt W$ with exceptional divisor
$\wt E\,=\,\sigma^{-1}(0)$.

We use the ``Theorem on formal functions'' \cite[Theorem~11.1]{Ha},
and the ``Grauert comparison theorem'' \cite[Theorem III.3.1]{BSt} for the analytic case.
Let ${\mathfrak m}_0 \,\subset \,\cO_{\wt W}$ be the maximal ideal sheaf for the
point $0\,\in\, \wt W$. Then the completion
$\left((R^1{\wt\sigma}_*\cO_{\wt S})_0\right)^{\!\boldsymbol{\vee}}$ of
$(R^1{\wt\sigma}_*\cO_{\wt S})_0$ in the $\mathfrak m_0$-adic topology is equal to
$$
\lim_{{\substack{\longleftarrow \\ k}}} H^1\left({\wt \sigma}^{-1}(0),\, \cO_{\wt S}/{\wt \sigma}^*(\mathfrak{m}^k_0)\right).
$$
We have the exact sequence
$$
0 \,\longrightarrow\,
\cO_E(k) \,\longrightarrow\, \cO_{\wt S}/{\wt \sigma}^*(\mathfrak m^{k+1}_0) \,
\longrightarrow\,\cO_{\wt S}/{\wt \sigma}^*(\mathfrak m^{k}_0) \,\longrightarrow\, 0
$$
of sheaves with support on
$$
{\wt \sigma}^{-1}(0)\,=\,\wt E\,\simeq\, \mathbb P^{n-1}
$$
so that the cohomology groups $H^q(\wt E,\, \cO_{\wt E}(k))$ vanish for all $k \geq 0$, and $q > 0$. In particular the maps
$$
H^1(\wt S,\, \cO_{\wt S}/{\wt \sigma}^*(\mathfrak m^{k+1}_0)) \,\longrightarrow \,
H^1(\wt S,\,\cO_{\wt S}/{\wt \sigma}^*(\mathfrak m^{k}_0))
$$
are isomorphisms for $k\geq 1$, and furthermore we have
$$
H^1(\wt S,\,\cO_{\wt S}/{\wt \sigma}^*(\mathfrak m_0))\,\simeq\, H^1(\mathbb P^{n-1},\,
\cO_{\mathbb P^{n-1}})\,=\,0\, .
$$
This shows
that $R^1{\wt \sigma}_* \cO_{\wt S}$ vanishes. This establishes the claim for blow-up at a point. 

Now consider the case where the center of the blow-up $\sigma$ is a smooth submanifold $A$ of
positive dimension. Since the claim is local with respect to the base space $W$, we may assume that $W$ is of the
form $A \times \wt W$, where both $A$ and $\wt W$ are small open subsets of complex number spaces,
e.g.\ polydisks. Denote by $\pi \,:\, W \,\longrightarrow\, \wt W$ the projection. We identify $A$ with
$A\times\{0\} =\pi^{-1}(0) \subset W$ as a submanifold.

Note that the blow-up
$$
\sigma\,:\,S \,\longrightarrow\, W
$$
of $W$ along $A$ is the fiber product
$\wt S \times_{\wt W} W \,\longrightarrow\, W$. The exceptional divisor $E$ of $\sigma$ can be identified with
$A\times \wt E$.

In the above argument we replace the maximal ideal sheaf $\mathfrak m_0$ by the vanishing ideal $\mathcal I_A$ of $A$.
Now $\sigma^*(\mathcal I^k_A)/\sigma^*(\mathcal I^{k+1}_A) \,\simeq\,\cO_{E}(k)$, and by
\cite[Theorem III.3.4]{BSt} we have
$$
R^1(\left.\sigma\right|_{E})_* \cO_{E} \,\simeq \,
\pi^* R^1(\left.\wt\sigma\right|_{\wt E})_*\cO_{\wt E}\,=\,0
$$
so that the
earlier argument can be applied. Hence the claim.

Now, let
$$
Z\,=\,Z_N \,\stackrel{\tau_N}{\longrightarrow} \,Z_{N-1} \,\stackrel{\tau_{N\!-\!1}}{\longrightarrow}\,
\cdots\, \stackrel{\tau_2}{\longrightarrow} \, Z_1\,\stackrel{\tau_1}{\longrightarrow}\,Z_0
\,=\, Y
$$
be the sequence of blow-ups that constitute $h\,:\,Z \,\longrightarrow\, Y$.
We have $\tau_{j*}\cO_{Z_j} \simeq \cO_{Z_{j-1}}$ and
$R^1\tau_{j\, *}\cO_{Z_{j}} = 0$ for $1\,\leq\, j\,\leq\, \tau_N$. Combining these
with \eqref{eq:leray} yields a canonical injective homomorphism
$$
H^1(Z_j, \,\cO_{Z_{j}})\,\longrightarrow\, H^1(Z_{j-1},\, \cO_{Z_{j-1}})
$$
that is an isomorphism for all $j\,=\,1,\, \cdots,\, N$. Hence, by naturality, the
homomorphism $\Theta_h\,:\,H^1(Y,\,\cO_Y) \,\longrightarrow\, H^1(Z,\,\cO_Z)$ is an isomorphism.
By our above remarks, this establishes the result.

\section*{Acknowledgements}

The second-named author thanks the Philipps-Universit\"at, where the work for this paper was carried
out, for its hospitality.

\end{document}